\documentclass{article}[12pt]
\usepackage{amsmath,amsfonts,amsthm}
\usepackage{amssymb}

\usepackage{epsfig}

\usepackage{subfigure}
\usepackage{graphics}
\usepackage{graphicx}
\usepackage{color}
\usepackage{soul}
\usepackage[makeroom]{cancel}

\usepackage{latexsym}
\usepackage{graphics}
\usepackage{amsmath}
\usepackage{amsthm}
\usepackage{xspace}
\usepackage{amssymb}
\usepackage{color}
\usepackage{cite}
\usepackage{latexsym}
\usepackage{graphics}
\usepackage{amsmath}
\usepackage{amsthm}
\usepackage{xspace}
\usepackage{amssymb}
\usepackage{epsfig}
\usepackage{comment}

\newcommand{\beql}[1]{\begin{equation}\label{#1}}
\newcommand{\eeql}{\end{equation}}
\newcommand{\eqn}[1]{(\ref{#1})}

\newcommand{\R}{\mathbb{R}}
\newcommand{\pr}{\mathbb{P}}
\newcommand{\E}{\mathbb{E}}

\newcommand{\cx}{{\cal X}}
\newcommand{\cm}{{\cal M}}

\newcommand{\cw}{{\cal W}}

\newcommand{\cc}{{\cal C}}

\newcommand{\bI}{{\bf I}}

\newtheorem{thm}{Theorem}
\newtheorem{lem}[thm]{Lemma}

\newtheorem{cor}[thm]{Corollary}

\newtheorem{definition}[thm]{Definition}

\theoremstyle{remark}
\newtheorem{remark}[thm]{Remark}

\begin{document}

\title{A particle system with mean-field interaction:\\
Large-scale limit of stationary distributions
%\\ (DRAFT) 
}

\author
{
Alexander L. Stolyar \\
ISE Department
and Coordinated Science Lab\\
University of Illinois at Urbana-Champaign\\
Urbana, IL 61801\\
\texttt{stolyar@illinois.edu}
}

\date{\today}

\maketitle

\begin{abstract}

We consider a system consisting of $n$ particles, moving forward in jumps on the real line.
System state is the empirical distribution of particle locations.
Each particle ``jumps forward''  at some time points, with the instantaneous rate of
jumps given by a decreasing function of the particle's location quantile within the current state (empirical distribution).
Previous work on this model established, under certain conditions, the convergence, as $n\to\infty$, of the system random
dynamics to that of a deterministic mean-field model (MFM), which is a solution to an integro-differential equation. Another line of previous work established the existence of MFMs that are traveling waves, as well as the attraction of MFM trajectories to traveling waves. The main results of this paper are: (a) We prove that, as $n\to\infty$, the stationary distributions of (re-centered) states concentrate on a (re-centered) traveling
wave; (b) We obtain a uniform across $n$ moment bound on the 
stationary distributions of (re-centered) states; (c) We prove a convergence-to-MFM result, which is substantially more general than that in previous work. Results (b) and (c) serve as ``ingredients'' of the proof of (a), but also are of independent interest.

\end{abstract}



%\newpage

\section{Introduction}

We consider a system consisting of $n$ particles, moving forward on the real line.
The particles move in jumps.
The system state at a given time is the current empirical distribution of particle locations.
Each particle gets ``urges to jump'' as an independent Poisson process of constant rate. 
However, a particle getting a jump urge actually jumps with the probability 
given by a decreasing function of the particle's location quantile within the current state (i.e., empirical distribution);
hence this a mean-field type of particles' interaction with each other. When a particle does jump, the jump size 
is independent, distributed as a random variable $Z>0$. We are interested in the system behavior when $n$ is large.

This model was introduced in \cite{GSS96, GMP97} as an idealized model
of distributed parallel simulation. In this case $n$ particles represent $n$ processors (``sites'') simulating different part of some large system, and a particle location is the current ``local simulation time'' of the corresponding processor. 
The following types of questions are of interest, as $n$ becomes large: how the local times of the processors progress over time;
 do local times ``stay closely together;''  does the evolution of the empirical distribution of local times becomes that of a traveling wave; etc. This model, and similar models, are motivated by other applications as well
 (including recent applications, such as blockchains), where, roughly speaking, a distributed 
 synchronization of a large number of sites is of interest; cf. \cite{Man-Sch-2005, Mal-Man-2006, Manita-2009, Manita-2006, Malyshkin-2006,Manita-2014,balazs-racz-toth-2014}, and references therein, for examples of synchronization models.

There are two lines of work on the particle system described above. The first is in paper \cite{GMP97},
where it was shown that, under certain additional conditions, as $n\to\infty$, the system random
dynamics converges to that of a deterministic mean-field model (MFM), which is a solution to an integro-differential equation.
There are several additional assumptions made in \cite{GMP97}, one of which is especially restrictive, in that the proof technique crucially relies on it -- a particle jump probability depends on the current locations of $K$ other particles chosen uniformly at random, where $K>0$ is a fixed parameter, same for all $n$. We will refer to this additional assumption as the {\em finite-dependence} assumption, and it substantially restricts the more general model of this paper. 

The second line of work is represented by papers \cite{GSS96, St2020-wave}, where formally defined mean-field models 
(solutions to an integro-differential equation) are studied. The results of \cite{GSS96} prove, in particular, that if 
an MFM that is a traveling wave exist, then ``typically'' the traveling wave is unique and MFM trajectories are attracted to
it as time increases; the question of a traveling wave existence under general assumptions was left open in \cite{GSS96}. 
Paper \cite{St2020-wave} proves the existence of a traveling wave, under very general assumptions. 
(We also note that, for some MFMs that are different from the one in this paper, 
the existence and explicit forms of the traveling waves are obtained in \cite{balazs-racz-toth-2014,Hongler2015,Hongler2019},
in some special cases of the jump size distribution.)

The combination of the results of \cite{GMP97,GSS96, St2020-wave} strongly suggests the following 
asymptotic property of the stationary distributions: as $n\to\infty$, the stationary distributions of (re-centered) states concentrate on a (re-centered) traveling wave. This property is stated as Conjecture 7.1 in \cite[Section 7]{St2020-wave}.

Our main results are:
\begin{description}
\item[(a)] We prove (Theorem~\ref{th-interchange}) that, as $n\to\infty$, the stationary distributions of (re-centered) states concentrate on a (re-centered) traveling wave. (This proves Conjecture 7.1 in \cite[Section 7]{St2020-wave}, under 
a slightly stronger assumption on the jump size, namely $\E Z^{2+\chi} < \infty$, $\chi >0$, as opposed to 
$\E Z^2 < \infty$.)

\item[(b)] As a key ``ingredient'' of the proof of (a), we obtain (Theorem~\ref{th-closeness}) a uniform across $n$ moment bound on the 
stationary distributions of (re-centered) states. This result is also of independent interest.

\item[(c)] We prove (Theorem~\ref{th-finite-interval}) the convergence-to-MFM result in our general setting, without the additional finite-dependence
assumption. (This substantially generalizes the result of \cite{GMP97}. The proof largely follows the approach used in 
\cite{balazs-racz-toth-2014}, for a different model. The approach is more generic than that in \cite{GMP97}; in particular,
it does not rely on the finite-dependence assumption.) This result is another ingredient of the proof of (a), but is also of independent interest.
\end{description}

The proof of (a) relies on (b) and (c), and on the results in \cite{GSS96, St2020-wave} on the existence/uniqueness of -- and attraction to -- traveling waves.

\subsection{Outline of the rest of the paper.}

In Section~\ref{sec-model} we formally define the model and informally state the main results.
Section~\ref{sec-notation} gives some basic notation, definitions and conventions used throughout the paper.
In Section~\ref{sec-formal-results} we state our main results, Theorems~\ref{th-interchange}, \ref{th-closeness},
and \ref{th-finite-interval} formally. Sections~\ref{sec-closeness-proof}, \ref{sec-finite-int-proof} and \ref{sec-proof-interchange}
contain the proofs of Theorems~\ref{th-closeness},
\ref{th-finite-interval} and \ref{th-interchange}, respectively.
A brief discussion of our results is in Section~\ref{sec-discuss}.

\section{The model and main results.}
\label{sec-model}

The particle system model is as follows. 
There are $n$ particles, moving in the positive direction (``right'') on the real axis $\R$. Each particle moves in jumps, as follows. For each particle there is an independent Poisson process of rate $\mu>0$ of ``jump urges.'' When a particle
gets an urge to jump, it actually jumps, to the right, with probability $\eta_n(\nu)$, where $\nu$ is its quantile in the current empirical distribution of the particles' locations;
that is, $\nu=\ell/n$ when the particle location is $\ell$-th from the left. With complementary probability $1-\eta_n(\nu)$ the particle does not jump.
To have the model well-defined, assume that quantile-ties between co-located particles are broken uniformly at random. 
We adopt the convention that function $\eta_n(\nu)$ is defined for a continuous argument $\nu \in [0,1]$, by assuming that 
it is constant in each interval $((\ell-1)/n, \ell/n]$ for $\ell=1,\ldots,n$ and $\eta_n(0)=1$. 
Assume that, for each $n$, function $\eta_n(\nu), ~0\le \nu \le 1,$ is non-increasing, and that, as $n\to\infty$, it uniformly converges 
to a continuous, strictly decreasing function $\eta(\nu), ~0\le \nu \le 1,$ 
with $\eta(0) = 1$, $\eta(1) = 0$. The jump sizes, when a particle does jump, are given by i.i.d. non-negative r.v. with CDF $J(y), y\ge 0$; we denote by $\bar J(y) = 1- J(y)$ the complementary CDF; a generic jump size is given by the r.v. $Z\ge 0$.
Without loss of generality, we can and will assume that $Z>0$, i.e. $J(0)=0$.
We will use notation 
\beql{eq-moment-finite222}
m^{(\ell)} \doteq \int_0^\infty y^\ell dJ(y) = \E Z^\ell, ~~\ell \ge 0,
\eeql
for the $\ell$-th moment of a jump size. (So, $m^{(1)}=\E Z$ is the mean jump size.) 

For some (not all) of our results we will need the following two additional conditions on the
 jump size distribution:
\beql{eq-moment-finite222}
m^{(2+\chi)} < \infty ~~\mbox{for some $\chi>0$}
\eeql
and
\beql{eq-positive-density}
\mbox{$J(\cdot)$ is absolutely continuous with density $J'(y)$, bounded away from $0$ on compact subsets of $\R_+$}.
\eeql

Note that, without loss of generality (WLOG), we can and do assume $\mu=1$; otherwise, we can achieve this condition by rescaling time. 
Also, if $m^{(1)} =\E Z < \infty$, we can and do assume WLOG that $m^{(1)} = 1$; otherwise, $m^{(1)} = 1$
is achieved by rescaling space.

Let $f^n(t) =(f_x^n(t), ~x\in \R)$ be the (random) empirical distribution of the 
particle locations at time $t$; namely, $f_x^n(t)$ is the fraction of particles located in $(-\infty,x]$ at time $t$.

As $n\to\infty$, it is very intuitive that $f^n_x(t)$ converges (in appropriate sense, under appropriate conditions) to a deterministic 
function $f_x(t)$, such that $f(t)=(f_x(t), ~x\in \R)$ is a distribution function for each $t$, and the following equation holds:
\beql{eq-dyn-trans-intro}
\frac{\partial}{\partial t} f_x(t) =
%- \mu \int_{-\infty}^x d_y f_y(t) \eta(f_y(t)) \bar J(x-y) =
-  \int_{-\infty}^x d_y f(y,t) \eta(f_y(t)) \bar J(x-y),
\eeql
where $d_y$ means the differential in $y$ (and recall that $\mu=1$ WLOG). 
We call a function $f_x(t)$ satisfying \eqn{eq-dyn-trans-intro} a {\em mean-field model.} 
(The formal meaning of \eqn{eq-dyn-trans-intro} and the definition of a mean-field model will be given 
later in Definition~\ref{def-mfm} in Section~\ref{sec-mfm}.)
The intuition for \eqn{eq-dyn-trans-intro} is as follows. 
For each $t$, the distribution $f(t)=(f_x(t), ~x\in \R)$ approximates the distribution of particles $f^n(t)$ when $n$ is large.
Since particles move right, $f_x(t)$ is non-increasing in $t$ for each $x$. So, the partial derivative 
$(\partial/\partial t) f_x(t)$ is non-positive and it should be equal to
the RHS of \eqn{eq-dyn-trans-intro}, which gives the instantaneous rate (scaled by $1/n$ and taken with minus sign) at which particles jump over point $x$ at time $t$.

It is known \cite{GSS96, St2020-wave} that, as long as $m^{(1)}<\infty$ (or, $m^{(1)}=1$ WLOG),
for any mean-field model,
 the speed at which the mean $\int x d_x f_x(t)$
of the distribution $f(t)$ moves right, is equal to $v = m^{(1)} \mu  \int_0^1 \eta(\nu) d\nu = \int_0^1 \eta(\nu) d\nu$. 
A distribution function $\phi=(\phi_x, ~x\in \R)$ is called a {\em traveling wave shape} (TWS) if $f_x(t) = \phi_{x-vt}$ is a mean-field model.
By substituting into \eqn{eq-dyn-trans-intro}, we see that any TWS $\phi$
must satisfy equation
\beql{eq-wave}
v  \phi'_x =
%\mu \int_{-\infty}^x \phi'_y \eta(\phi_y) \bar J(x-y) dy = 
\int_{-\infty}^x \phi'_y \eta(\phi_y) \bar J(x-y) dy.
\eeql
It is known \cite{St2020-wave} that a TWS $\phi$ exists as long as $m^{(2)} < \infty$ (no other assumption 
on the jump size distribution is required), and it is unique (up to a shift) if, in addition, \eqn{eq-positive-density} holds.

We now informally state the main results of this paper. (The formal results will be given later, in Section~\ref{sec-formal-results},
after introducing more notation.)
Let $\mathring f^n(t)$ denote the distribution $f^n(t)$, re-centered so that the distribution mean is at $0$.

%\begin{prop}[The formal result is Theorem~\ref{th-interchange}]
%\label{prop-interchange}
{\bf Informal statement of Theorem~\ref{th-interchange} (in Section~\ref{sec-formal-results}).}
{\em Assume \eqn{eq-moment-finite222} and \eqn{eq-positive-density}. Then, as $n\to\infty$, 
the stationary distribution of the process $\mathring f^n(\cdot)$ converges to the (Dirac) distribution
concentrated on the unique TWS $\phi$, centered so that its mean is at $0$. Moreover, $\phi$
has finite absolute $(1+\chi)$-th moment: $\int_{-\infty}^\infty |y|^{1+\chi} d\phi_y < \infty$.}
%\end{prop}

%\begin{prop}[The formal result is Theorem~\ref{th-closeness}]
%\label{prop-closeness}
{\bf Informal statement of Theorem~\ref{th-closeness} (in Section~\ref{sec-formal-results}).}
{\em Assume \eqn{eq-moment-finite222} and \eqn{eq-positive-density}.
Then, for all sufficiently large $n$, the Markov process $\mathring f^n(\cdot)$ is stable
(positive Harris recurrent) and its stationary distribution is such that 
the expected absolute $(1+\chi)$-th moment of $\mathring f^n(\cdot)$ is bounded uniformly in $n$:
$$ %\beql{eq-tight-uniform}
\E \int_{-\infty}^\infty |y|^{1+\chi} d_y \mathring f^n_y(t) \le \bar C.
%\E \Phi_{1+\chi}(\mathring f^n(\infty)) \le \bar C.
$$
} %\eeql
%\end{prop}

%\begin{prop}[The formal result is Theorem~\ref{th-finite-interval}]
%\label{prop-finite-interval}
{\bf Informal statement of Theorem~\ref{th-finite-interval} (in Section~\ref{sec-formal-results}).}
{\em Assume $m^{(1)} <\infty$ (or, $m^{(1)}=1$ WLOG) and \eqn{eq-positive-density}. 
Suppose that, as $n\to\infty$, the initial conditions $f^n(0)$ converge to a deterministic proper
distribution $f(0)$. (Nothing else about $f(0)$ is assumed, not even the existence of a finite mean.)
Then the process $f^n(\cdot)$ converges to the unique mean-field model $f(\cdot)$
with initial condition $f(0)$.}
%\end{prop}

\section{Basic notation}
\label{sec-notation}

The set of real numbers is denoted by $\R$, and is viewed as the usual Euclidean space.
As a measurable space, $\R$ is endowed with Borel $\sigma$-algebra. % $\mathcal B(\R)$. 
For scalar functions $h(x)$ of a real $x$: $\|h\|_1 = \int_x |h(x)|dx$ is $L_1$-norm; 
$h(x)$ is called $c$-Lipschitz if it is Lipschitz with constant $c \ge 0$. 
Let $\cc_b$ be the set of continuous bounded functions on $\R$, which are constant outside a closed interval (one constant value to the ``left'' of it, and possibly another constant value to the ``right'' of it.)

For functions $h(x)$ of a real $x$: 
$h(x+)$ and $h(x-)$ are the right and left limits; 
a function $h(x)$ is RCLL if it is right-continuous and has left limits at each $x$. 

A function $h$ of $x$ may be written as either $h(x)$ or $h_x$.
Notation $d_x h(x,t)$ for a multivariate function $h(x,t)$, where $x\in \R$, 
 means the differential in $x$.
 
Denote by $\cm$ the set of scalar RCLL non-decreasing
 functions $f=(f(x), ~x\in \R)$, which are (proper) probability distribution functions,
i.e., such that $f(x) \in [0,1]$, $\lim_{x\downarrow -\infty} f(x) = 0$ and $\lim_{x\uparrow \infty} f(x) = 1$.
For elements $f \in \cm$ we use the terms {\em distribution function} and {\em distribution} interchangeably.
Space $\cm$ is endowed with Levy-Prohorov metric (cf. \cite{Ethier_Kurtz}) and the corresponding topology
of weak convergence (which is equivalent to the convergence at every point of continuity of the limit);
the weak convergence in $\cm$ is denoted $\stackrel{w}{\rightarrow}$.
Note that, for  $f,\phi \in \cm$, the $L_1$-norm of their difference, $\|f-\phi\|_1$, is equal to the Wasserstein 
$W_1$-distance between the corresponding two distributions. 
The inverse ($\nu$-th quantile) 
of $f \in \cm$ is $f^{-1}(\nu) \doteq \inf\{y~|~f(y) \ge \nu\}$, $\nu\in [0,1]$; $\gamma^{-1}(1) =\infty$ when 
$f(y) < 1$ for all $y$.

Unless explicitly specified otherwise, we use the following conventions regarding random elements and random processes.
A measurable space is considered equipped with a Borel $\sigma$-algebra, induced by the 
metric which is clear from the context. A random process $Y(t), ~t\ge 0,$ always takes values in a complete separable metric space (clear
from the context), and has RCLL sample paths.
For a random process $Y(t), ~t\ge 0,$ we denote by $Y(\infty)$ the random value of $Y(t)$ in a stationary regime (which will be clear from the context). Symbol $\Rightarrow$ signifies convergence of random elements in distribution; $\stackrel{\pr}{\longrightarrow}$ means convergence in probability.
 {\em W.p.1} or {\em a.s.} means {\em with probability one.}
{\em I.i.d.} means {\em independent identically distributed.}
 For a condition/event $A$, $\bI\{A\}=1$ if $A$ holds, and $\bI\{A\}=0$ otherwise.

Space $D([0,\infty), \R)$ [resp. $D([0,\infty), \cm)$] is the Skorohod space of RCLL functions on $[0,\infty)$ taking values in 
$\R$ [resp. $\cm$], with the corresponding Skorohod ($J_1$) metric and topology (cf. \cite{Ethier_Kurtz});
$\stackrel{J_1}{\rightarrow}$ denotes the convergence in these spaces.

For a distribution $f \in \cm$ and scalar function $h(x), x \in \R$, 
$f h \doteq \int_{\R} h(x) df_x$. 

For scalar functions $h(x), x \in \cx$, with some domain $\cx$, $\|h\| = \sup_{x\in \cx} |h(x)|$ is the sup-norm.
When $G_k, G$ are operators mapping the space of such functions into itself, 
$\lim G_k h = Gh$ and $G_k h \to Gh$  mean the uniform convergence: 
$\|G_k h - Gh\| \to 0$.

Suppose we have a Markov process with state space $\cx$ and transition function $P^t(x,H),$ $t\ge 0$. 
A measurable set $X \subseteq \cx$ is called {\em small} (cf. \cite{Bramson-book}),
if there exist constants $T>0$ and $\delta>0$, and probability distribution $\alpha(\cdot)$ on $\cx$,
such that for any $x \in X$ and any measurable $H\subseteq \cx$, $P^T(x,H) \ge \delta \alpha(H)$.
$P^t$, as an operator, is $P^t h(x)  \doteq \int_y P^t(x,dy) h(y)$, where $h$ is a scalar function with domain $\cx$;
 $I=P^0$ is the identity operator.
The process (infinitesimal) generator $B$ is 
$$
Bh \doteq \lim_{t\downarrow 0} (1/t) [P^t - I] h.
$$
Function $h$ is within the domain of the generator $B$ if $Bh$ is well-defined. We say that a Markov process is {\em stable}
if it is {\em positive Harris recurrent} (cf. \cite{Bramson-book}); if it is, it has unique stationary distribution.

 For real numbers $a$ and $b$
 we use notations: $\mbox{sign}(a) \doteq \bI(a\ge 0) - \bI(a < 0)$, 
$a\wedge b \doteq \min\{a,b\}$.
RHS and LHS mean right-hand side and left-hand side, respectively; WLOG means {\em without loss of generality}.
Abbreviation {\em w.r.t.} means {\em with respect to};
{\em a.e.} means {\em almost everywhere w.r.t. Lebesgue measure}.

\section{Formal statements of main results}
\label{sec-formal-results}

\subsection{Uniform moment bound for stationary distributions. Limit of stationary distributions.} 

For an element $f\in \cm$, denote by $\bar f$ the mean of the corresponding distribution,
$$
\bar f = \int_{-\infty}^\infty x df_x = \int_0^\infty (1-f_x) dx - \int_{-\infty}^0 f_x dx,
$$
with the usual convention that the mean is well-defined and finite when both integrals in the RHS are finite.
Denote
$$
\mathring \cm = \{f \in \cm ~|~ \bar f =0\}.
$$
When $\bar f$ is finite, denote by $\mathring f = (\mathring f_x,~x\in \R) \in \mathring \cm$ the centered version of $f$, namely
$$
\mathring f_x = f_{x+\bar f}, ~x\in \R.
$$
If we denote by $\cm^{(n)} \subset \cm$ the state space of the process $f^n(\cdot)$, then 
$\mathring \cm^{(n)} = \cm^{(n)} \cap \mathring \cm$ is the state space of $\mathring f^n(\cdot)$. 
If we use notation
$$
\Phi_\ell(f) = \int_{-\infty}^\infty |x|^\ell d f_x,
$$
for the $\ell$-th absolute moment of $f$, then we obviously have
$$
\Phi_\ell(f) < \infty, ~~~\forall f \in \cm^{(n)}, ~\forall n, ~\forall \ell \ge 0.
$$

\begin{thm}
\label{th-interchange}
Suppose conditions \eqn{eq-moment-finite222} and \eqn{eq-positive-density} hold.
Then, as $n\to\infty$, 
\beql{eq-main-conv}
\mathring f^n(\infty) \Rightarrow \phi,
\eeql
where $\phi$ is the unique TWS with $\bar \phi = 0$. Moreover,
$\phi$ is such that 
\beql{eq-tws-moment}
\Phi_{1+\chi}(\phi)
% = \int_{-\infty}^\infty |y|^{1+\chi} d\phi_y 
< \infty,
\eeql
and a stronger form of convergence \eqn{eq-main-conv} holds: 
\beql{eq-main-conv-stronger}
\| \mathring f^n(\infty) - \phi \|_1 \Rightarrow 0.
\eeql
\end{thm}

Proof of Theorem~\ref{th-interchange} is in Section ~\ref{sec-proof-interchange}.

\begin{thm}
\label{th-closeness}
Suppose conditions \eqn{eq-moment-finite222} and \eqn{eq-positive-density} hold.
Then there exist $\bar C>0$ and $\bar n$ such that, for all $n \ge \bar n$ 
the Markov process $\mathring f^n(\cdot)$ is stable and we have
\beql{eq-tight-uniform}
\E \Phi_{1+\chi}(\mathring f^n(\infty)) \le \bar C.
\eeql
\end{thm}

Proof of Theorem~\ref{th-closeness} is in Section~\ref{sec-closeness-proof}.

\subsection{Transient behavior: convergence to a mean-field model
}

Essentially all our asymptotic results on the transient behavior of the processes 
are for the {\em non-centered} processes $f^n(\cdot)$. 
The result for the centered processes (Theorem~\ref{th-finite-interval}(ii)) is obtained essentially as a corollary. 

Denote by $L^{(n)}$ the generator of the process
$f^n(\cdot)$. 
For any $h\in\cc_b$, function $f^n h$ of $f^n$ is within the domain of $L^{(n)}$ (where we use the fact each function in $\cc_b$ is constant outside a closed interval), and
$$
L^{(n)} [f^n h] = 
\int_0^1 d\nu [f^n]^{-1}(\nu) \eta_n([f^n]^{-1}(\nu)) \E [h([f^n]^{-1}(\nu)+Z) - h([f^n]^{-1}(\nu))]
$$
where the expectation is over the distribution of the random jump size $Z$.  
We also formally define the ``limit'' of $L^{(n)}$ as
$$
L [f h] = 
\int_0^1 d\nu f^{-1}(\nu) \eta(f^{-1}(\nu)) \E [h(f^{-1}(\nu)+Z) - h(f^{-1}(\nu))]
$$

\begin{thm}
\label{th-finite-interval}
Suppose $f^n(0) \stackrel{w}{\rightarrow} f(0)$, where $\{f^n(0)\}$ is deterministic sequence of $f^n(0) \in \cm^{(n)}$,
and $f(0) \in \cm$. (Note that we do {\em not} assume that $f(0)$ has well-defined mean $\bar f(0)$.)
Assume that conditions $m^{(1)}=\E Z<\infty$ (or, $m^{(1)} = \E Z =1$ WLOG)
and  \eqn{eq-positive-density} hold. Then we have: 

(i) $f^n(\cdot) \Rightarrow  f(\cdot)$ in $D([0,\infty), \cm)$, where $f(\cdot) \in D([0,\infty), \cm)$ is deterministic, uniquely determined by $f(0)$. Moreover, $f(\cdot)$ is a continuous element of $D([0,\infty), \cm)$, which satisfies 
\beql{eq-pde-oper}
f(t) h - f(0) h - \int_0^t L f(s) h ds =0, ~~\forall h\in \cc_b, ~\forall t\ge 0.
\eeql
The dependence of $f(\cdot)$ (as an element of space $D([0,\infty), \cm)$ with $J_1$-convergence topology) on $f(0)$ (as an element of $\cm$ with weak convergence topology) is continuous.

(ii) If, in addition, $\bar f(0) =0$, then $\bar f(t) =vt$ and consequently $\mathring f(\cdot)$ is a continuous element of $D([0,\infty), \cm)$, uniquely determined by $\mathring f(0)=f(0) \in \mathring \cm$, with $\mathring f(t) \in \mathring \cm$
for all $t\ge 0$. The dependence of $\mathring f(\cdot)$  on $\mathring f(0)$ is continuous.
And if, in addition, $\bar f^n(0) =0$ for all $n$, then 
$\mathring f^n(\cdot) \Rightarrow \mathring f(\cdot)$ in $D([0,\infty), \cm)$.

\end{thm}

Proof of Theorem~\ref{th-finite-interval} is in Section~\ref{sec-finite-int-proof}. Also in Section~\ref{sec-finite-int-proof}
we show (Theorem~\ref{lem-fde-oper-equiv}) that solutions to \eqn{eq-pde-oper} are exactly the mean-field models 
(Definition~\ref{def-mfm}).
We note that many of the supplementary results in Section~\ref{sec-finite-int-proof}, 
which may be of independent interest, require
assumptions on the jump size distribution that are {\em much weaker} than conditions $m^{(1)} = \E Z < \infty$ and \eqn{eq-positive-density}. In particular, some of those result assume nothing about the jump size distribution besides it being a proper distribution.

\section{Proof of Theorem~\ref{th-closeness}}
\label{sec-closeness-proof}

\subsection{Equivalent view of process $\mathring f^n(\cdot)$.}
\label{sec-equiv-view}

State $\mathring f^n(t)$ can be equivalently described as
$w^n(t)=(w_1(t), w_2(t), \ldots, w_n(t)) \in \R^n$, where 
$w_1(t), w_2(t), \ldots, w_n(t)$ are the locations of the $n$ particles 
w.r.t. the mean $\bar f^n(t)$,
listed in a non-decreasing order.
 (So, the average $(1/n)\sum_i w_n(t)=0$ at all times.)
From now on for each $\mathring f^n \in \mathring\cm^{(n)}$ we will consider the corresponding vector 
$w^n=(w_1, w_2, \ldots, w_n)$, and vice versa. Any function of 
$\mathring f^n \in \mathring\cm^{(n)}$ may be expressed via the corresponding $w^n$, and vice versa. 
In particular, 
$$
\Phi_\ell(\mathring f^n(t)) = 
\frac{1}{n} \sum_{i=1}^n |w_i(t)|^\ell.
$$
Note that the topology on $\mathring\cm^{(n)}$, induced by the (weak convergence) topology on $\cm$,
is equivalent to the usual topology of component-wise convergence of the corresponding vectors $w^n$.

The evolution of $w^n(t)$ is as follows. 
Between the times of the jump urges, $w^n(t)$ remains constant.
At a time $t$ of a jump urge, the following occurs. Let $\kappa_i(t)$ be the actual jumps size of particle $i$,
in the system without re-centering, upon this urge; $\kappa_i(t)\ge 0$ and can be non-zero for at most one particle.
Then, in the re-centered system, particle $i$ jump size (i.e., the increment of $w_i(t)$) at $t$ 
is $\zeta_i=\kappa_i(t) - \sum_s \kappa_s(t)/n$ (which may be positive or negative). After the jumps (if any) at $t$ occur,
the particles indices $i$ are changed, if necessary, to keep $w_i(t)$ non-decreasing in $i$.

\subsection{Informal intuition for the proof.}
\label{sec-tight-intuition}

The proof of stability, in Subsection~\ref{sec-stabil}, uses fluid limit technique and is fairly straightforward. Let us discuss the intuition for the proof of the bound \eqn{eq-tight-uniform} in Subsection~\ref{sec-moment-uniform}.

At a very high level, the bound \eqn{eq-tight-uniform} is due to the fundamental property of the system, which can be called the ``egalitarian trend:'' in re-centered system, the particles at high quantiles (large $i$) will have a negative drift, while particles at low quantiles (small $i$) will have positive drift, thus preventing 
the centered empirical distribution $\mathring f^n$ from ``spreading out.''

To obtain the bound on the expected $(1+\chi)$-th moment of $\mathring f^n$, we need the finite $(2+\chi)$-th moment on a jump size. Informally speaking, we use
$\Phi_{2+\chi}(\mathring f^n)$ as a Lyapunov function. If $B^{(n)}$ is the generator of process $\mathring f^n(\cdot)$, then, ``generally speaking,'' we have $\E B^{(n)} \Phi_{2+\chi}(\mathring f^n(\infty))=0$; in other words,
the expected drift of $\Phi_{2+\chi}(\mathring f^n(t))$ in steady-state is $0$. 
Function
$$
G_{1+\chi}^{(n)} (\mathring f^n)
= (2+\chi)  \sum_i \mbox{sign}(w_i) |w_i|^{1+\chi} \E  \zeta_i^{(\mathring f^n)}, ~~ \mathring f^n \in \mathring\cm^{(n)},
$$
where $\zeta_i^{(\mathring f^n)}$ is the random jump size of particle $i$ (in re-centered system) upon a jump urge when the state is 
$\mathring f^n$, can be thought of as the ``first-order approximation of the generator $B^{(n)}$, applied to function $\Phi_{2+\chi}(\mathring f^n)$;'' note that the derivative $|w_i^{2+\chi}|' = (2+\chi) \mbox{sign}(w_i) w_i^{1+\chi}$. We can show that,
when $\Phi_{1+\chi}(\mathring f^n)$ is large, 
$$ %\beql{eq-drift-33}
%G_{1+\chi}^{(n)} (\mathring f^n) = 
G_{1+\chi}^{(n)} (\mathring f^n) \le -\epsilon \Phi_{1+\chi}(\mathring f^n)
$$ %\eeql
for some constant $\epsilon>0$; this is where we use the egalitarian trend property, which ensures that $\E  \zeta_i^{(\mathring f^n)}$ is negative [resp., positive]
for particles at high [resp., low] quantiles. From here we can obtain, informally speaking,
$$
G_{1+\chi}^{(n)} (\mathring f^n) \le -\epsilon \Phi_{1+\chi}(\mathring f^n) +K,
$$
which holds for some $K>0$ and all $\mathring f^n$. Taking into account the fact that $G_{1+\chi}^{(n)} (\mathring f^n)$
is not $B^{(n)}$, but only its first-order approximation, and doing the corresponding estimates, we obtain,
informally speaking,
$$
B^{(n)} \Phi_{2+\chi}(\mathring f^n) \le -(\epsilon/2) \Phi_{1+\chi}(\mathring f^n) +K.
$$
Taking expectation w.r.t. $\mathring f^n(\infty)$, we obtain, informally speaking,
$$
0 \le -(\epsilon/2) \E \Phi_{1+\chi}(\mathring f^n(\infty)) +K,
$$
which yields \eqn{eq-tight-uniform}.

In the actual proof, instead of $\Phi_{2+\chi}(\mathring f^n)$
we use its truncated version $\Phi_{2+\chi}^{(C_1)}(\mathring f^n) = \Phi_{2+\chi}(\mathring f^n) \wedge C_1$
as the Lyapunov function, because the latter is certainly within the domain of generator $B^{(n)}$. And then let $C_1\uparrow \infty$.

We note that, this ``program'' for proving a property of the type of \eqn{eq-tight-uniform} is likely applicable to 
other models having the egalitarian trend property, while the technical details may differ.

\subsection{Stability.}
\label{sec-stabil}

The stability (positive Harris recurrence) is easily established using the fluid limit technique (\cite{RS92, Dai95, St95, Bramson-book}).
We note that the stability proof only uses conditions \eqn{eq-positive-density} and $m^{(1)} < \infty$ 
(as opposed to stronger condition \eqn{eq-moment-finite222}).

To prove stability we will use the equivalent representation $w^n(\cdot)$ of the process $\mathring f^n(\cdot)$,
given in Subsection~\ref{sec-equiv-view}.
The state space of $w^n(\cdot)$ is 
$$\mathring\cw^{(n)} = \{ w^n=(w_1, \ldots, w_n) \in \R^n ~|~ (1/n)\sum_i w_i=0\},$$
i.e. the set of those vectors in $\R^n$, corresponding to $\mathring f^n \in \mathring\cm^{(n)}$.
The norm of a state $w^n=(w_1, \ldots, w_n) \in \mathring\cw^{(n)} $  is $\|w^n\|= \max_i |w_i|$.
Using \eqn{eq-positive-density}, it is straightforward to see that, for any fixed $a >0$,
the closed set $\mathring\cw^{(n)}(a)=\{w^n \in \mathring\cw^{(n)}  ~|~\|w^n\| \le a\}$ is small (see the definition 
in Section~\ref{sec-notation}).

Pick any two numbers, $0< \nu_1< \nu_2 < 1$, and  choose $\bar n$ large enough so that for any $n \ge \bar n$,
$\eta_n(\nu_1) - \eta_n(\nu_2) > (\eta(\nu_1) - \eta(\nu_2))/2$. (All we need here is that the jump probabilities of the ``left-most''
and ``right-most'' particle are separated by a positive constant.) Consider any fixed $n \ge \bar n$.

Consider a sequence of versions of the process $w^n(\cdot)$, namely 
processes $w^{n,k}(\cdot)$ with increasing norm of the initial state, 
$\|w^{n,k}(0)\| = c_k \uparrow \infty$, $k\to\infty$, with $w^{n,k}(0)/c_k \to w(0) = (w_1(0), \ldots, w_n(0))$.
Given that $\mathring\cw^{(n)}(a)$ is a closed small set for any $a>0$,
to establish stability it suffices to show that, for some fixed $T>0$, $\|w^{n,k}(T)\|/c_k \Rightarrow 0$.
This in turn follows from the fact that the limit (in appropriate sense) $w(\cdot)=(w_1(\cdot), \ldots, w_n(\cdot))$
of the sequence of processes $w^{n,k}(c_k t)/c_k, ~t\ge 0,$ has trajectories such that 
$w(t) \in \mathring\cw^{(n)}$ for all $t\ge 0$,
$\|w(0)\|=1$ and
$(d/dt) \max w_i(t) \le - \epsilon <0$ as long as $\max w_i(t) >0$; therefore, $\|w(t)\|=0$ for all $t \ge 1/\epsilon$.
We omit further details, which are straightforward. $\Box$

\subsection{Proof of \eqn{eq-tight-uniform}} 
\label{sec-moment-uniform}

At some level, this proof is similar to the proof of an analogous result in \cite{St2021-coc} for a different particle system. However,
the difference of our model from that in \cite{St2021-coc} is substantial, so we give full details of the proof for our model.

Consider the following function 
$$
G_{1+\chi}^{(n)} (\mathring f^n)
= (2+\chi)  \sum_i \mbox{sign}(w_i) |w_i|^{1+\chi} \E  \zeta_i^{(\mathring f^n)}, ~~ \mathring f^n \in \mathring\cm^{(n)},
$$
where $\zeta_i^{(\mathring f^n)}$ is the random jump size (which can have any sign) of particle $i$ upon a jump urge when the state is 
$\mathring f^n$. (The sizes $\zeta_i^{(\mathring f^n)}$ are dependent across $i$, of course.) 

As we will see, function $G_{1+\chi}^{(n)} (\mathring f^n)$ can be thought of as the ``first-order approximation of the generator $B^{(n)}$ of process $\mathring f^n(\cdot)$, applied to function $\Phi_{2+\chi}(\mathring f^n)$;'' but we do not even claim that $\Phi_{2+\chi}(\mathring f^n)$ is within the generator $B^{(n)}$ domain.
Note that, for each $n$, $G_{1+\chi}^{(n)} (\mathring f^n)$ is continuous in $\mathring f^n \in \mathring\cm^{(n)}$.

Note that each 
$\E  \zeta_i^{(\mathring f^n)}$ is the quantity of the order $O(1/n)$, which motivates the definition
\beql{eq-bar-zeta-def}
\bar \zeta_i^{(\mathring f^n)} \doteq n \E \zeta_i^{(\mathring f^n)}.
\eeql
We observe that, if particle $i$ location $w_i$ is the $\ell$-th from the left, and it is not co-located with any other particle,
then
$$
\bar \zeta_i^{(\mathring f^n)} = \eta_n(\ell/n) -v_n, 
$$
where 
\beql{eq-vn-def}
v_n \doteq \sum_\ell \eta_n(\ell/n) \equiv \int_0^1 \eta_n(\nu) d\nu
\eeql
is the average drift of the mean of the {\em non-centered} particle system. (Clearly, $\lim_n v_n = v$.) In the more general
case, when exactly $k$ particles are co-located particles -- namely $\ell$-th, $(\ell+1)$-th, ..., $(\ell+k-1)$-th
left-most particles are co-located -- and particle $i$ is one of them, we have
\beql{eq-zeta-vn}
\bar \zeta_i^{(\mathring f^n)} = \frac{\eta_n(\ell/n) + \ldots + \eta_n((\ell+k-1)/n)}{k} -v_n.
\eeql
We define the function $\bar\zeta^{(\mathring f^n)}(x), ~x\in \R,$ as follows: 
$\bar \zeta^{(\mathring f^n)}(x) = \bar \zeta^{(\mathring f^n)}_i$, where $i$ is the particle whose location $w_i$ is the closest to $x$ on the left; we also adopt a convention that, if $w_i$ is the location 
of the left-most particle, then $\bar\zeta^{(\mathring f^n)}(x) = \bar \zeta_i^{(\mathring f^n)}$ for all $x< w_i$.
Clearly, function $\bar\zeta^{(\mathring f^n)}(x)$ is a piece-wise constant non-increasing function.

We can write:
\beql{eq-gen1-33}
G_{1+\chi}^{(n)} (\mathring f^n)
= (2+\chi)  \int_{-\infty}^{\infty} \bar\zeta^{(\mathring f^n)}(x) \mbox{sign}(x) |x|^{1+\chi} d\mathring f^n_x,
\eeql
or, by ``integrating over the values $\nu$ of  $\mathring f^n_x$'' and using \eqn{eq-zeta-vn},
\beql{eq-gen1-alt-33}
G_{1+\chi}^{(n)} (\mathring f^n)
= (2+\chi)  \int_0^1 d\nu [\eta_n(\nu) - v_n] \mbox{sign}([\mathring f^n]^{-1}(\nu)) |[\mathring f^n]^{-1}(\nu)|^{1+\chi} \le 0.
\eeql
The inequality in \eqn{eq-gen1-alt-33} follows by the following argument. Denote 
$$
\nu_0^n \doteq \inf\{\nu~|~\eta_n(\nu)  < v_n\},
$$
and observe that, as $n\to\infty$, 
$$
\lim_{n\to\infty} \nu_0^n = \nu_0, ~\mbox{where} ~ \eta(\nu_0)= v.
$$
To be specific, consider the case when $[\mathring f^n]^{-1}(\nu_0^n) \le 0$. (The case $[\mathring f^n]^{-1}(\nu_0^n) \ge 0$
is treated analogously.) We have $\nu_*^n \doteq \mathring f^n_0 \ge \nu_0^n$. Note that
$$
\int_0^{\nu^n_0} [\eta_n(\nu) - v_n] d\nu  = - \int_{\nu^n_0}^1 [\eta_n(\nu) - v_n] d\nu.
$$
Then we have
$$
\int_{\nu^n_*}^1 d\nu [\eta_n(\nu) - v_n] \mbox{sign}([\mathring f^n]^{-1}(\nu)) |[\mathring f^n]^{-1}(\nu)|^{1+\chi} =
$$
$$
\int_{\nu^n_*}^1 d\nu [\eta_n(\nu) - v_n] |[\mathring f^n]^{-1}(\nu)|^{1+\chi} \le 0,
$$
and
$$
\int_0^{\nu^n_*} d\nu [\eta_n(\nu) - v_n] \mbox{sign}([\mathring f^n]^{-1}(\nu)) |[\mathring f^n]^{-1}(\nu)|^{1+\chi} =
$$
$$
- \int_0^{\nu^n_*} d\nu [\eta_n(\nu) - v_n]  |[\mathring f^n]^{-1}(\nu)|^{1+\chi} =
$$
$$
- \int_0^{\nu^n_0} d\nu [\eta_n(\nu) - v_n] |[\mathring f^n]^{-1}(\nu)|^{1+\chi} 
- \int_{\nu^n_0}^{\nu^n_*} d\nu [\eta_n(\nu) - v_n] |[\mathring f^n]^{-1}(\nu)|^{1+\chi} \le 0,
$$
where the last inequality holds because $|[\mathring f^n]^{-1}(\nu)|^{1+\chi}$ is non-increasing in $[0,\nu^n_*]$.
Thus, \eqn{eq-gen1-alt-33} is proved.

Next, we claim the following property: there exists a sufficiently large $C>0$ and some $\epsilon > 0$, such that,
uniformly in all sufficiently large $n$ and all $\mathring f^n \in \mathring\cm^{(n)}$  with $\Phi_{1+\chi}(\mathring f^n) \ge C$,
\beql{eq-drift-33}
G_{1+\chi}^{(n)} (\mathring f^n) \le -\epsilon \Phi_{1+\chi}(\mathring f^n).
\eeql
The proof of \eqn{eq-drift-33} is given in Section~\ref{sec-eq-drift-33}.

From \eqn{eq-drift-33} and \eqn{eq-gen1-alt-33} we obtain that, uniformly in all sufficiently large $n$,
\beql{eq-drift2-33}
G_{1+\chi}^{(n)}(\mathring f^n) \le -\epsilon \Phi_{1+\chi}(\mathring f^n) +\epsilon C.
\eeql

Denote by $\Phi_{2+\chi}^{(C_1)}(\mathring f^n) = \Phi_{2+\chi}(\mathring f^n) \wedge C_1$ the function $\Phi_{2+\chi}$ truncated at level $C_1>0$.
Given that this is a continuous bounded function of $\mathring f^n \in \mathring \cm^{(n)}$, and the process of jump urges is Poisson,
it is not hard to see that $\Phi_{2+\chi}^{(C_1)}(\mathring f^n)$ is within the domain 
of the generator $B^{(n)}$ of process $\mathring f^n(\cdot)$.

Next, we claim the following fact:
there exist $C_2>0$ such that for any fixed $C_1>0$, 
uniformly in all large $n$ and $\mathring f^n$ such that $\Phi_{2+\chi}(\mathring f^n) \le C_1$, we have
\beql{eq-generator-bound-33}
B^{(n)} \Phi_{2+\chi}^{(C_1)}(\mathring f^n) \le G_{1+\chi}^{(n)}(\mathring f^n) + C_2 \Phi_\chi(\mathring f^n) + C_2,
\eeql
and then
\beql{eq-generator-bound5-33-int}
B^{(n)} \Phi_{2+\chi}^{(C_1)}(\mathring f^n) \le  -\epsilon \Phi_{1+\chi}(\mathring f^n) + C_2 \Phi_\chi(\mathring f^n) + C_3,
\eeql
with $C_3 = \epsilon C + C_2$.
The proof of \eqn{eq-generator-bound-33} is given in Section~\ref{sec-generator-bound-33}.

From \eqn{eq-generator-bound5-33-int} and inequality $\Phi_\chi(\mathring f^n) \le [\Phi_{1+\chi}(\mathring f^n)]^{\chi/(1+\chi)} \le C_4 + (\epsilon/(2 C_2)) \Phi_{1+\chi}(\mathring f^n)$, which holds for a sufficiently large fixed $C_4>0$,
we obtain
\beql{eq-generator-bound5-33}
B^{(n)} \Phi_{2+\chi}^{(C_1)}(\mathring f^n) \le  -(\epsilon/2) \Phi_{1+\chi}(\mathring f^n) + C_5,
\eeql
where $C_5 = C_4 C_2 + C_3$.

Bound \eqn{eq-generator-bound5-33} in turn implies that
 for any fixed $C_1>0$, 
\beql{eq-generator-bound2-33}
B^{(n)} \Phi_{2+\chi}^{(C_1)}(\mathring f^n)  \le [-(\epsilon/2) \Phi_{1+\chi}(\mathring f^n)  + C_5] \bI \{\Phi_{2+\chi}(\mathring f^n) \le C_1\},
\eeql
because, obviously, $B^{(n)} \Phi_{2+\chi}^{(C_1)}(\mathring f^n) \le 0$ when $\Phi_{2+\chi}(\mathring f^n) > C_1$.

Recalling that $\mathring f^n(\infty)$ is the random value of $\mathring f^n(t)$ in the stationary regime, 
we have for all large $n$:
$$
0 = \E B^{(n)}  \Phi_{2+\chi}^{(C_1)}(\mathring f^n(\infty)) \le 
\E \left[ (-(\epsilon/2) \Phi_{1+\chi}(\mathring f^n(\infty)) + C_5) \bI \{\Phi_{2+\chi}(\mathring f^n(\infty)) \le C_1\} \right]
$$
$$
\le -(\epsilon/2) \E [\Phi_{1+\chi}(\mathring f^n(\infty))  \bI \{\Phi_{2+\chi}(\mathring f^n(\infty)) \le C_1\}] 
+ C_5
$$
and then
$$
\E [\Phi_{1+\chi}(\mathring f^n(\infty))  \bI \{\Phi_{2+\chi}(\mathring f^n(\infty))  \le C_1\}] \le 2 C_5/\epsilon.
$$
Letting $C_1 \uparrow \infty$, we finally obtain that
$$
\E \Phi_{1+\chi}(\mathring f^n(\infty)) \le 2 C_5/\epsilon
$$
for all sufficiently large $n$, and then 
$$
\E \Phi_{1+\chi}(\mathring f^n(\infty)) \le \bar C
$$
holds for all $n$ for some large $\bar C>0$.
$\Box$

\subsection{Proof of \eqn{eq-drift-33}.}
\label{sec-eq-drift-33}

The definition of $\bar \zeta_i=\bar \zeta_i^{(\mathring f^n)}$ in \eqn{eq-bar-zeta-def} can be interpreted as follows:  
$\bar \zeta_i$ is the expected jump size of particle $i$,  
{\em conditioned on this particle receiving the jump urge}, and then centered by the expected 
jump size of any particle upon a jump urge in the system.

The proof is by contradiction. Suppose property \eqn{eq-drift-33} does not hold. Then, we can and do choose a subsequence of $n\to\infty$, and corresponding $\mathring f^n$, so that 
along this subsequence $\Phi_{1+\chi}(\mathring f^n) \uparrow \infty$ and 
\beql{eq-contr-33}
G_{1+\chi}^{(n)} (\mathring f^n) / \Phi_{1+\chi}(\mathring f^n) \to 0.
\eeql

Consider separately two cases (a) $\liminf_n \Phi_1(\mathring f^n) = c <\infty$ and 
(b) $\lim_n \Phi_1(\mathring f^n) = \infty$.

Case (a).

Consider a subsequence of $\mathring f^n$ such that $\lim_n \Phi_1(\mathring f^n) = c <\infty$ and, moreover, 
$\mathring f^n \stackrel{w}{\rightarrow} f$, where $f$ is a proper distribution.

For a fixed $0<\delta < 1/2$ denote
$$
\Phi_{1+\chi}^{(\delta)}(\mathring f^n) = 
\int_{[0,\delta] \cup [1-\delta,1]} \mbox{sign}([\mathring f^n]^{-1}(\nu))  |[\mathring f^n]^{-1}(\nu)|^{1+\chi} d\nu.
$$
Using the facts that $f$ is proper and $\Phi_{1+\chi}(f^n) \to \infty$, we easily see that, for any fixed $\delta>0$,
$$
\lim \Phi_{1+\chi}^{(\delta)}(\mathring f^n) / \Phi_{1+\chi}(\mathring f^n) =1.
$$
Pick a sufficiently small $\delta>0$ such that, for some $\delta_1>0$ and all large $n$,  
$\eta_n(\delta) -v_n > \delta_1$ and $\eta_n(1-\delta) -v_n < -\delta_1$. Then, we have
$$
\liminf |G_{1+\chi}^{(n)} (\mathring f^n)|/\Phi_{1+\chi}(\mathring f^n) \ge (2+\chi)\delta_1,
$$
which contradicts \eqn{eq-contr-33}.

Case (b). 

Denote $c_n=\Phi_1(\mathring f^n) \uparrow \infty$, and consider the sequence of rescaled versions of $\mathring f^n$, namely
$$
\tilde f^n_x = \mathring f^n_{c_n x}, ~~w\in \R.
$$
Note that 
\beql{eq-same-ratios}
G_{1+\chi}^{(n)} (\mathring f^n)/\Phi_{1+\chi} (\mathring f^n) = G_{1+\chi}^{(n)} (\tilde f^n)/\Phi_{1+\chi} (\tilde f^n).
\eeql
Consider two subcases: (b.1) $\lim_n \Phi_{1+\chi}(\tilde f^n) = \infty$ and (b.2) $\liminf_n \Phi_{1+\chi}(\tilde f^n) = c <\infty$. In the subcase (b.1), on account of \eqn{eq-same-ratios}, we can obtain a contradiction in the same way as in the case (a). 

So, the remaining case to consider is (b.2). 

For an $f \in \cm$, let us formally define a ``limiting version'' of the functional $G_{1+\chi}^{(n)}(\mathring f^n)$ defined in 
the LHS of the inequality in  \eqn{eq-gen1-alt-33}:
\beql{eq-gen1-alt-33-lim}
G_{1+\chi} (f) 
= (2+\chi)  \int_0^1 d\nu [\eta(\nu) - v] \mbox{sign}(f^{-1}(\nu)) |f^{-1}(\nu))|^{1+\chi}.
\eeql

Note that $\Phi_{1+\chi}(\tilde f^n) \ge 1$, so that $0<c<\infty$.
Consider a subsequence of $\tilde f^n$ such that  $\lim_n \Phi_{1+\chi}(\tilde f^n) = c$ and
$\tilde f^n \stackrel{w}{\rightarrow} f \in \cm$. Distribution $f$ cannot be concentrated at a single point. (Otherwise, since $\Phi_1(\tilde f^n)=1$ for all $n$, $\Phi_{1+\chi}(\tilde f^n)$ could not remain bounded.) 
Therefore, $|G_{1+\chi}(f)| > 0$, and then
$$
\liminf |G_{1+\chi}^{(n)} (\tilde f^n)| \ge |G_{1+\chi}(f)| > 0,
$$
and then $\liminf |G_{1+\chi}^{(n)} (\tilde f^n)|/\Phi_{1+\chi}(\tilde f^n) >0$, which contradicts \eqn{eq-contr-33}.
$\Box$

\subsection{Proof of \eqn{eq-generator-bound-33}.}
\label{sec-generator-bound-33}

We will use the following inequality, which holds, for some constant $C_6>0$, for any numbers $y$ and $\delta$:
\beql{eq-number-ineq}
|y+\delta|^{2+\chi} -|y|^{2+\chi} 
\le (2+\chi) \mbox{sign}(y) |y|^{1+\chi} \delta
+ C_6 |y|^\chi \delta^2 + C_6 |\delta|^{2+\chi}.
\eeql
Indeed, 
$$
|y+\delta|^{2+\chi} -|y|^{2+\chi} \le (2+\chi) \mbox{sign}(y) |y|^{1+\chi} \delta
+ (1/2)(2+\chi)(1+\chi)  |\tilde y|^\chi \delta^2,
$$
for some $\tilde y \in [y-|\delta|,y+|\delta|]$. Using 
$|\tilde y|^\chi \le (|y|+|\delta|)^\chi \le (2|y|)^\chi + (2|\delta|)^\chi$,
we obtain \eqn{eq-number-ineq}.

Consider a fixed state $\mathring f^n$ and consider the expected increment $\Delta$ of $\Phi_{2+\chi}(\mathring f^n)$ upon a jump urge, which occurs in this state. 
Then, using \eqn{eq-number-ineq},
$$
\Delta= \E \frac{1}{n} \sum_i |w_i+\zeta_i|^{2+\chi} - \frac{1}{n} \sum_i |w_i|^{2+\chi} \le 
\E \frac{1}{n} \sum_i [(2+\chi)~\mbox{sign}(w_i) |w_i|^{1+\chi} \zeta_i + C_6 |w_i|^\chi \zeta_i^2 + C_6 |\zeta_i|^{2+\chi}],
$$
where expectation $\E$ is with respect to the uniform selection of the particle receiving the jump urge, the random event of it actually jumping, and the randomness of the jump size (if it occurs).

For $\zeta_i$ we have:
$$
\zeta_i = \kappa_i - \frac{1}{n} \sum_s \kappa_s, 
$$
where $\kappa_i=\kappa_i(\mathring f^n)$ the (random) jump size of particle $i$.
Note that $\E \zeta_i$ is the quantity of the order $O(1/n)$, since $\sum_s \kappa_s$ is of order $O(1)$
and $\E \kappa_i $ is of order $O(1/n)$ (because $1/n$ is the probability of particle $i$ receiving the jump urge). 
Therefore, $\bar \zeta_i = n \E \zeta_i = O(1)$,
and we can write:
$$
\E \sum_i \mbox{sign}(w_i) |w_i|^{1+\chi} \zeta_i =  \frac{1}{n} \sum_i \mbox{sign}(w_i) |w_i|^{1+\chi} \bar \zeta_i.
$$
Next,
$$
\E (\zeta_i/2)^{2+\chi} \le (1/2)  \E \kappa_i^{2+\chi}  +
(1/2)  \frac{1}{n^{2+\chi}} \E (\sum_s \kappa_s)^{2+\chi}
$$
and therefore
\beql{eq-zeta-moment}
\E \zeta_i^{2+\chi} \le C_7/n
\eeql
for some $C_7>0$
because $\E (\sum_s \kappa_s)^{2+\chi}$ is upper bounded by the $(2+\chi)$-th moment 
of a jump size and $\E \kappa_i^{2+\chi}$ is upper bounded by $1/n$ times the $(2+\chi)$-th moment 
of a jump size. Note that \eqn{eq-zeta-moment} holds for $\chi=0$ as well, with possibly different $C_7$.
Therefore, by choosing $C_7$ sufficiently large, we have both \eqn{eq-zeta-moment} and
\beql{eq-zeta-moment2}
\E \zeta_i^{2} \le C_7/n.
\eeql

Assembling these bounds, we obtain
\beql{eq-n-times-drift-33}
n\Delta 
\le (2+\chi) \frac{1}{n} \sum_i \mbox{sign}(w_i) |w_i|^{1+\chi} \bar \zeta_i
+ C_6 C_7 \frac{1}{n} \sum_i |w_i|^\chi + 
C_6 C_7 = G_{1+\chi}^{(n)}(\mathring f^n) + C_2 \Phi_\chi(\mathring f^n) + C_2,
\eeql
where, recall, $\Delta$ is the expected increment  of $\Phi_{2+\chi}$ upon a jump urge occurring in a fixed state $\mathring f^n$,
and $C_2 = C_6 C_7$.

Now consider the value of the generator $B^{(n)} \Phi_{2+\chi}^{(C_1)}(\mathring f^n)$ at point $\mathring f^n$ such that $\Phi_{2+\chi}(\mathring f^n) \le C_1$.
 For that, consider 
the expected increment of $\Phi_{2+\chi}^{(C_1)}(\mathring f^n(t))$ over a small interval $[0,t/n]$, with $\mathring f^n(0)=\mathring f^n$.
First, note that, as $t\downarrow 0$, the contribution into this expected increment of the event that more than one jump urge occurs, is $o(t)$. (Because jump urges follow a Poisson process of rate $n$, and $\Phi_{2+\chi}^{(C_1)}$ is bounded.) With probability
$t + o(t)$ there will be exactly one jump urge in $[0,t/n]$, which therefore occurs into the state $\mathring f^n$ (such that $\Phi_{2+\chi}(\mathring f^n) \le C_1$); then, the expected increment of $\Phi_{2+\chi}^{(C_1)}$ will not exceed that
of $\Phi_{2+\chi}$.
Using these observations and the estimate \eqn{eq-n-times-drift-33}, we obtain \eqn{eq-generator-bound-33}.
We omit the remaining straightforward $\epsilon/\delta$
formalities. 
$\Box$

\section{Proof of Theorem~\ref{th-finite-interval}}
\label{sec-finite-int-proof}

This proof largely follows the approach used in \cite{balazs-racz-toth-2014}, for a different model. 
Unlike in \cite{balazs-racz-toth-2014}, we work with Levy-Prohorov metric 
(inducing the weak convergence topology) on $\cm$, as opposed to
the stronger Wasserstein $W_1$-metric. This, in fact, simplifies some parts of the proof in our case;
we will point out those parts as they appear.
However, some other parts of our proof of Theorem~\ref{th-finite-interval} are completely different from (or not present in)
the development in \cite{balazs-racz-toth-2014}. They are: Section~\ref{sec-mfm} and Theorem~\ref{lem-fde-oper-equiv},
which establish the equivalence between solutions to \eqn{eq-pde-oper} and mean-field models;
Theorem~\ref{lem-fde-unique}, which establishes the uniqueness of a mean-field model and its continuous dependence
on the initial state.

We note that many of the supplementary results in this section, which may be of independent interest, require
assumptions on the jump size distribution that are {\em much weaker} than conditions $m^{(1)} = \E Z < \infty$ and \eqn{eq-positive-density}. In particular, some of the result assume nothing about the jump size distribution besides it being a proper distribution. We will emphasize such weaker assumptions, where applicable, in the results' statements. 

\subsection{$C$-relative compactness of the processes.}

A sequence of random processes with sample path in the Skorohod space $D([0,\infty), \R)$ 
[resp., $D([0,\infty), \cm)$] is called {\em $C$-relatively compact} (see \cite{perkins-2002,balazs-racz-toth-2014})
if it is: (a) relatively compact, i.e. its any subsequence has a further subsequence converging in distribution to some limiting process; and (b) any such limiting process has continuous sample paths, a.s.

\begin{thm}
\label{th-6.9-analog}
Suppose $f^n(0) \stackrel{w}{\rightarrow} f(0)$, where $\{f^n(0)\}$ is deterministic sequence of $f^n(0) \in \cm^{(n)}$,
and $f(0) \in \cm$. 
Then for any $h \in \cc_b$ the sequence of processes $\{f^n(\cdot) h\}$ is
$C$-relatively compact in the Skorohod space $D([0,\infty), \R)$. 
(Note that we do {\em not} assume that $f(0)$ has well-defined mean $\bar f(0)$, or \eqn{eq-positive-density},
or \eqn{eq-moment-finite222}, or even $m^{(1)} = \E Z < \infty$. Jump size distribution only needs to be proper.)
\end{thm}

{\em Proof.} This result is analogous to Theorem 6.9 in \cite{balazs-racz-toth-2014}.
Note that, although the model in \cite{balazs-racz-toth-2014} is different from ours, the only property
that is used in the proof of Theorem 6.9 in \cite{balazs-racz-toth-2014} is that the rate of jumps of each particle is upper 
bounded by a finite constant $a$ at all times. The later property, obviously, holds for our model as well, with $a=\mu=1$.
Therefore, the proof of Theorem 6.9 in \cite{balazs-racz-toth-2014} applies essentially verbatim, with the following 
adjustments. 

What in the proof of Theorem 6.9 in \cite{balazs-racz-toth-2014} are $f, \mu_n, a$, in our notation are $h,f^n, \mu=1$, respectively.  In the proof, $x_i(t)$ denote the locations of the particles, uniquely determined by the process state 
at time $t$, and vice versa. The notation $\tilde x_i(t)$ is used for the locations of particles in an artificial system,
with the same initial particle locations $\tilde x_i(0)=x_i(0)$, but such that each particle jumps every time it gets a jump urge;
the artificial system is coupled to the original one in the natural way, so that the corresponding particles have common processes of jump urges and common jump sizes (if the particle in the original system happen to jump at a jump urge). Clearly,
with this coupling, $x_i(t) \le \tilde x_i(t)$ and $x_i(t) - x_i(s) \le \tilde x_i(t) - \tilde x_i(s)$ at all times $s \le t$.
Finally, in our case, $I_n(t) = f^n(t) h = \int h(x) d_x f_x^n(t) = (1/n)\sum_i h(x_i(t))$.
$\Box$

\begin{thm}
\label{th-6.10-analog}
Suppose $f^n(0) \stackrel{w}{\rightarrow} f(0)$, where $\{f^n(0)\}$ is deterministic sequence of $f^n(0) \in \cm^{(n)}$,
and $f(0) \in \cm$. 
 Then the sequence of processes $\{f^n(\cdot)\}$ is
$C$-relatively compact in the Skorohod space $D([0,\infty), \cm)$.
(Note that we do {\em not} assume that $f(0)$ has well-defined mean $\bar f(0)$, or \eqn{eq-positive-density},
or \eqn{eq-moment-finite222}, or even $m^{(1)} = \E Z < \infty$. Jump size distribution only needs to be proper.)
\end{thm}

{\em Proof.} 
This result is analogous to Corollary 6.10 in \cite{balazs-racz-toth-2014}, with essentially same proof.
In fact, in our case the proof is simpler. We need to verify conditions (i) and (ii) of Theorem II.4.1 in  \cite{perkins-2002},
which in our case take the following form.

(i) For any $T>0$ and $\epsilon>0$, there exists  $K>0$ such that
$$
\sup_n \pr \left(\sup_{t\le T}    [f^n_{-K}(t) + 1- f^n_{K}(t)] >\epsilon  \right) < \epsilon.
$$
(ii) For any $h \in \cc_b$ the sequence of processes $\{f^n(\cdot) h\}$ is
$C$-relatively compact in the Skorohod space $D([0,\infty), \R)$. 
(Note that the class of functions $h \in \cc_b$
is {\em separating}, which means that a 
probability measure $g$ is uniquely determined by
the values of $g h$ for $h\in \cc_b$.)

Condition (ii) is verified by Theorem~\ref{th-6.9-analog}. The verification of condition (i) repeats  the proof of 
Corollary 6.10 in \cite{balazs-racz-toth-2014}, essentially verbatim, up to and including the display where Markov inequality is used for the first time. At that point it remains to observe that the probability in the RHS of the display can be made arbitrarily small by making $K$ sufficiently large. 
%Specifically, the adjustments of the proof of Corollary 6.10 in \cite{balazs-racz-toth-2014} are as follows.
The measure $\mu_n(t,\cdot)$ in the proof of Corollary 6.10 in \cite{balazs-racz-toth-2014} is in our notation the measure (distribution) $f^n(t)$ on $\R$;
%given by the system state at time $t$, i.e. by the distribution function $(f_x^n(t), x \in \R)$ in our notation. 
the particle locations $x_i(t)$ and $\tilde x_i(t)$ in the coupled original and artificial systems,
are as described above in the proof of our Theorem~\ref{th-6.10-analog}.
$\Box$

\subsection{Trajectories of a limit satisfy \eqn{eq-pde-oper}
}

For trajectories $f^n(\cdot) \in D([0,\infty), \cm)$ with $f^n(t) \in \cm^{(n)}$ for all $t\ge 0$, 
let us define the following functional for each $h\in \cc_b$ and $t\ge 0$:
$$
A^n_{t,h}(f^n(\cdot)) \doteq f^n(t) h - f^n(0) h - \int_0^t L^{(n)} f^n(s) h ds.
$$
We will also formally define a ``limit version'' of $A^n_{t,h}$ for trajectories $f(\cdot) \in D([0,\infty), \cm)$, 
for each $h\in \cc_b$ and $t\ge 0$:
$$
A_{t,h}(f(\cdot)) \doteq f(t) h - f(0) h - \int_0^t L f(s) h ds.
$$

\begin{thm}
\label{th-6.11-analog}
Suppose $f^n(0) \stackrel{w}{\rightarrow} f(0)$, where $\{f^n(0)\}$ is deterministic sequence of $f^n(0) \in \cm^{(n)}$,
and $f(0) \in \cm$.  
Then the sequence of processes $\{f^n(\cdot)\}$ is such that,
for every $t \ge 0$ and any $h \in \cc_b$,
$$
\sup_{0 \le s \le t} |A^n_{s,h}(f^n(\cdot))| \Rightarrow 0, ~\mbox{as}~ n\to\infty.
$$
(Note that we do {\em not} assume that $f(0)$ has well-defined mean $\bar f(0)$, or \eqn{eq-positive-density},
or \eqn{eq-moment-finite222}, or even $m^{(1)} = \E Z < \infty$. Jump size distribution only needs to be proper.)
\end{thm}

{\em Proof.} The proof repeats 
the proof of  Theorem 6.11 in \cite{balazs-racz-toth-2014}, in which we replace: $L$ by $L^{(n)}$; $f$ by $h$; $\mu_n$
by $f^n$; $a$ by $\mu=1$. In particular, in our case, $I_n(t) = f^n(t)h = \int h(x) d_x f_x^n(t) = (1/n)\sum_i h(x_i(t))$, so that
$L I_n(t)$ is replaced by
$$
L^{(n)} I_n(t) = L^{(n)} f^n(t)h = \int_0^1 \E[h(f_{\nu}^{-1}(t)+Z)-h(f_{\nu}^{-1}(t))]  \eta_n(\nu) d\nu,
$$
where the expectation in the integrand is over a random jump size $Z$. The martingale $M_n(t), ~t\ge 0,$ in our case is
$$
M_n(t) = A^n_{t,h}(f^n(\cdot)),
$$
so that $L M^2_n(t)$ in our case is $L^{(n)} M^2_n(t)$. The last line of the last display of the proof in \cite{balazs-racz-toth-2014}
can be removed, and the final estimate $L^{(n)} M^2_n(t) \le 4\|h\|^2/n$ for $h \in \cc_b$ can be observed without that
line (because the jump urge rate of each particle is $\mu=1$).

Finally, note that in our theorem we only need to consider $h \in \cc_b$. We do {\em not} need 
to consider the identity test function $h(x)=x$, and that is why in Theorem~\ref{th-6.11-analog} we 
do {\em not} need condition $\E Z^2 < \infty$ -- it suffices that $Z$ has a proper distribution. $\Box$

\begin{thm}
\label{th-6.12-analog}
Suppose $f^n(0) \stackrel{w}{\rightarrow} f(0)$, where $\{f^n(0)\}$ is deterministic sequence of $f^n(0) \in \cm^{(n)}$,
and $f(0) \in \cm$.
Suppose, the sequence of processes $\{f^n(\cdot)\}$ is such that,
$f^n(\cdot) \Rightarrow f(\cdot)$ as $n\to\infty$, in $D([0,\infty), \cm)$.
Then, for every $t \ge 0$ and any $h \in \cc_b$,
$$
A^n_{t,h}(f^n(\cdot)) \Rightarrow A_{t,h}(f(\cdot)), ~\mbox{as}~ n\to\infty.
$$
(Note that we do {\em not} assume that $f(0)$ has well-defined mean $\bar f(0)$, or \eqn{eq-positive-density},
or \eqn{eq-moment-finite222}, or even $m^{(1)} = \E Z < \infty$. Jump size distribution only needs to be proper.)
\end{thm}

{\em Proof.} The statement of this theorem is analogous to that of Theorem 6.12 in \cite{balazs-racz-toth-2014}.
However, the proof in our case is simpler, and is as follows.
We know that, by Theorem~\ref{th-6.10-analog}, a.s. the limiting process $f(\cdot)$ has continuous
trajectories in $D([0,\infty), \cm)$.
We can use Skorohod representation to construct all processes on a common probability space so that, w.p.1, 
$f^n(\cdot) \stackrel{J_1}{\rightarrow} f(\cdot)$ as $n\to\infty$; moreover, by the continuity of $f(\cdot)$,
we see that $f^n(t) \to f(t)$ uniformly on compact sets of $t$; we also see that, w.p.1,
$f_x(t)$ is non-increasing in $t$ (because so are $f^n_x(t)$).
Then, it is easy to see (using, in particular, the facts that
convergence $\eta_n(\cdot) \to \eta(\cdot)$ is uniform, and $\eta(\cdot)$ is strictly decreasing continuous) that
$A^n_{t,h}(f^n(\cdot)) \to A_{t,h}(f(\cdot))$ w.p.1 $\Box$

\begin{cor}
\label{cor-6.13-analog}
Suppose $f^n(0) \stackrel{w}{\rightarrow} f(0)$, where $\{f^n(0)\}$ is deterministic sequence of $f^n(0) \in \cm^{(n)}$,
and $f(0) \in \cm$. 
Then the sequence of processes $\{f^n(\cdot)\}$ is such that
its any distributional limit $f(\cdot)$ in $D([0,\infty), \cm)$ is such that, a.s., $f(\cdot)$ is continuous 
and it satisfies $A_{t,h}(f(\cdot))=0$ (i.e., \eqn{eq-pde-oper}) 
for every $t \ge 0$ and any $h \in \cc_b$.
(Note that we do {\em not} assume that $f(0)$ has well-defined mean $\bar f(0)$, or \eqn{eq-positive-density},
or \eqn{eq-moment-finite222}, or even $m^{(1)} = \E Z < \infty$. Jump size distribution only needs to be proper.)
\end{cor}

{\em Proof.} By Theorem~\ref{th-6.10-analog}, there exists a subsequence of $\{f^n(\cdot)\}$, which converges in distribution 
to a process $f(\cdot)$ with a.s. continuous trajectories. By Theorems~\ref{th-6.11-analog} and \ref{th-6.12-analog},
this process must satisfy $A_{t,h}(f(\cdot))=0$
for every $t \ge 0$ and any $h \in \cc_b$.
$\Box$

\subsection{Equivalent characterization of solutions to \eqn{eq-pde-oper} as mean-field models.}
\label{sec-mfm}

\begin{definition}
\label{def-mfm}
A function $f_x(t), ~x\in \R, ~t\in \R_+,$ will be called a {\em mean-field model} if it satisfies the following conditions.\\
 (a) For any $t$, $f(t) = (f_x(t), x\in \R) \in \cm$.\\
(b) For any $x$, $f_x(t)$ is non-increasing $c$-Lipschitz in $t$, with constant $c$ independent of $x$.\\
(c) For any $x$, for any $t$ where the partial derivative
$(\partial/\partial t) f_x(t)$ exists (which is almost all $t$ w.r.t. Lebesgue measure, by the Lipschitz property), equation 
\beql{eq-dyn-trans}
\frac{\partial}{\partial t} f_x(t) = 
- \int_0^1 d\nu f^{-1}(\nu) \eta(f^{-1}(\nu)) 
\bI\{f^{-1}(\nu) \le x\} \bar J(x-f^{-1}(\nu))
\eeql
holds. 
\end{definition}

Note that, by the change of variable $y = f^{-1}(\nu)$ in the integral in the RHS of \eqn{eq-dyn-trans},
equation \eqn{eq-dyn-trans} can be equivalently written as
\beql{eq-dyn-trans-equiv}
\frac{\partial}{\partial t} f_x(t) = 
- \int_{-\infty}^x d_y f_y(t) \bar \eta(y,f(t)) \bar J(x-y),
\eeql
where we use notations
$$
    \bar \eta(y,f(t)) \doteq \left\{\begin{array}{ll}
        \eta(f_y(t)), & \text{when $f_u(t)$ is continuous at point $u=y$}\\
        (\nu_2-\nu_1)^{-1} \int_{\nu_1}^{\nu_2} \eta(\nu) d\nu, & \text{otherwise},
        \end{array}\right. 
$$
where $\nu_2=f_y(t), \nu_1=f_{y-}(t)$. 
Equation \eqn{eq-dyn-trans-equiv} (or \eqn{eq-dyn-trans})
is a more general form of \eqn{eq-dyn-trans-intro}, allowing $f_x(t)$ to be RCLL in $x$, rather than continuous.
If $f_u(t)$ is continuous at $u=y$, then  $\bar \eta(y,f(t))=\eta(f_y(t))$;
if $f_u(t)$ has a jump at $u=y$ then $\bar \eta(y,f(t))$ is
$\eta(\nu)$ averaged over $\nu\in [f_{y-}(t),f_y(t)]$. In paper \cite{St2020-wave}, the equation in form 
 \eqn{eq-dyn-trans-equiv} is used to define a mean-field model.

\begin{thm}
\label{lem-fde-oper-equiv}
For any initial condition $f(0) \in \cm$, $f(\cdot)$ satisfies \eqn{eq-pde-oper} if and only if it is a mean-field model. (Note that we do {\em not} assume that $f(0)$ has well-defined mean $\bar f(0)$, or \eqn{eq-positive-density},
or \eqn{eq-moment-finite222}, or even $m^{(1)} = \E Z < \infty$. Jump size distribution only needs to be proper.)
\end{thm}

{\em Proof.} 
`Only if.' Let $h=(h(u), ~u\in \R)$ be the left-continuous step-function $h(u) = \bI(u \le x)$,
jumping from 1 to 0 at point $x$.
Let $h_\epsilon$, $\epsilon>0$, be a continuous approximation of $h$, which 
is linearly decreasing from $1$ to $0$ in $[x,x+\epsilon]$. We see that 
$$
L [f(t) h_\epsilon] \to L [f(t) h], ~~\forall t.
$$
Indeed, $|L [f(t) h_\epsilon] - L [f(t) h]|$ is upper bounded by the probability that a random jump of size $Z$ of a particle 
 randomly located according to $f(t)$, is such that the jump either originates or lands in $(x,x+\epsilon)$;
 this probability vanishes as $\epsilon\downarrow 0$. Also, since both $L [f(t) h_\epsilon]$ and $L [f(t) h_\epsilon]$
 are within $[-1,0]$, for all $t$ and $\epsilon$, we have a universal bound $|L [f(t) h_\epsilon] - L [f(t) h]| \le 1$. 
 Then, for any fixed $t$,
by taking the $\epsilon\downarrow 0$ limit in $f(t) h_\epsilon - f(0) h_\epsilon - \int_0^t L [f(s) h_\epsilon] ds =0$, we obtain
$$
f(t) h  - f(0) h  - \int_0^t L [f(s) h] ds =0.
$$
This means that $f(t) h = f_x(t)$ is absolutely continuous in $t$, with the derivative equal to 
$(\partial/\partial t) f_x(t)=L [f(t) h]$ a.e. in $t$. This, in particular, implies that, for any fixed $y$, the $f_y(t)-f_{y-}(t)$
is continuous in $t$ (in fact, Lipschitz); this, in turn, means that, possible discontinuity points $y$ of $f_y(t)$
``cannot move''  in time $t$. We can then conclude that, for any $x$, the derivative
$(\partial/\partial t) f_x(t)=L [f(t) h]$ is in fact continuous in $t$.
Therefore, $(\partial/\partial t) f_x(t)=L [f(t) h]$ at every $t$. It remains to observe that
$L [f(t) h]$ is exactly the RHS of \eqn{eq-dyn-trans}.

`If.'  The definition of a mean-field model $f(\cdot)$ implies that \eqn{eq-pde-oper} holds for the defined above step-function $h$ for any $x$.
Then, we have \eqn{eq-pde-oper} for any $h$,
which is piece-wise constant with finite number of pieces; and the set of such functions $h$ is tight within the space
of test functions $h \in \cc_b$, equipped with uniform metric. Then \eqn{eq-pde-oper} holds for any $h \in \cc_b$.
$\Box$

\subsection{Uniqueness of solution to \eqn{eq-pde-oper} and continuity in initial state. Proof of Theorem~\ref{th-finite-interval}(i).}
\label{sec-uniqueness}

\begin{thm}
\label{lem-fde-unique}
Assume that $m^{(1)}=\E Z<\infty$ (or, $m^{(1)}=\E Z=1$ WLOG)
and condition \eqn{eq-positive-density} holds.
Then, for any initial condition $f(0) \in \cm$, a solution $f(\cdot)$ of \eqn{eq-pde-oper} (i.e., a mean-field model) is unique.
(Note that we do {\em not} assume that $f(0)$ has well-defined mean $\bar f(0)$, 
or \eqn{eq-moment-finite222}.)
\end{thm}

{\em Proof.}  
We know that solutions to \eqn{eq-pde-oper} are mean-field models.
Papers \cite{GSS96,St2020-wave} study the properties of mean-field models.
It is easy to check that the proofs of all results in Sections 4.4 and 4.5 of \cite{GSS96}, 
including Theorem 2, stating that the Wasserstein 
$W_1$-distance (i.e., the $L_1$-norm of the difference) between any two mean-field models $f^{(1)}(\cdot)$ and $f^{(2)}(\cdot)$ 
is non-increasing, {\em never use} the fact that the means $\bar f^{(1)}(0)$ and $\bar f^{(2)}(0)$ 
are well-defined and equal to $0$. Those proofs only use the fact that $f^{(1)}(0),f^{(2)}(0) \in \cm$ and 
$$
\int_{-\infty}^{\infty} [f_w^{(1)}(0)-f_w^{(2)}(0)] dw = 0.
$$
(The conditions $\E Z < \infty$ and 
\eqn{eq-positive-density}
{\em are used} there. Technically, those proofs in \cite{GSS96} assume
 that the jumps size distribution $J(\cdot)$ is exponential -- 
but only the property \eqn{eq-positive-density} is actually used.)

For any $f(0) \in \cm$ a mean-field model $f(\cdot)$ is such that (see \cite{St2020-wave})  
\beql{eq-speed-conserv-222}
\int_{-\infty}^{\infty} [f_w(0)-f_w(t)] dw = vt, ~~\forall t \ge 0.
\eeql
Now, if $f^{(1)}(\cdot)$ and $f^{(2)}(\cdot)$ are two different solutions with the same initial condition,
$f^{(1)}(0)=f^{(2)}(0)=f(0)$, then, from \eqn{eq-speed-conserv-222},
$$
\int_{-\infty}^{\infty} [f_w^{(1)}(t)-f_w^{(2)}(t)] dw = 0, ~~\forall t \ge 0.
$$
Therefore, the Wasserstein $W_1$-distance between $f^{(1)}(t)$ and $f^{(2)}(t)$ cannot increase, which implies the uniqueness.
$\Box$

{\em Proof of Theorem~\ref{th-finite-interval}(i).} We have established that the family of
distributions of the processes is $C$-tight, 
any subsequential distributional limit is concentrated on solutions $f(\cdot)$ to \eqn{eq-pde-oper} (i.e., mean-field models),
and the solution $f(\cdot)$ with initial condition $f(0)$ is unique. This implies the convergence $f^n(\cdot) \Rightarrow  f(\cdot)$.

The continuity of $f(\cdot)$ in $f(0)$ is then a consequence of convergence and uniqueness. Indeed, consider 
a sequence of initial states $f^{(k)}(0) \in \cm$, converging to some $f(0) \in \cm$; namely $f^{(k)}(0) \stackrel{w}{\rightarrow} f(0)$
as $k\to \infty$. Fix $\epsilon>0$ and any $\delta>0$.
 We can choose an increasing sequence $n=n(k) \to \infty$, and corresponding $f^{(k),n}(0) \in \cm^{(n)}$,
such that $f^{(k),n}(0) \stackrel{w}{\rightarrow} f(0)$ as $k\to\infty$, and for each $k$
the process $f^{(k),n}(\cdot)$ is within distance $\epsilon$ from the deterministic trajectory $f^{(k)}(\cdot)$
with probability at least $1-\delta$. 
Then, for all large $k$, $f^{(k),n}(\cdot)$ is within distance $\epsilon$ from both $f^{(k)}(\cdot)$ and $f(\cdot)$
with probability at least $1-2\delta$. Since this is true for any $\epsilon$ and $\delta$, 
the only option is that $f^{(k)}(\cdot) \stackrel{J_1}{\rightarrow} f(\cdot)$ .
$\Box$

\subsection{Corollaries for the centered processes: Proof of Theorem~\ref{th-finite-interval}(ii)}

Suppose $\bar f(0) =0$. 
Then, from \eqn{eq-speed-conserv-222} we have $\bar f(t) =vt$.
Then the continuity and uniqueness of $\mathring f(\cdot)$, as well as the continuity of its dependence on 
$\mathring f(0)$ follow from the corresponding properties of non-centered $f(\cdot)$ in Theorem~\ref{th-finite-interval}(i).

Suppose, in addition, that $\bar f^n(0) =0$ for all $n$.
Since $\eta_n$ converges to $\eta$, we easily see directly that $\bar f^n(t) \stackrel{\pr}{\rightarrow} vt$ for any $t$.
Combining these observations with Theorem~\ref{th-finite-interval}(i), we obtain 
the convergence $\mathring f^n(\cdot) \Rightarrow \mathring f(\cdot)$;
for example, we can use Skorohod representation to obtain a.s. convergence to the limit. 
$\Box$

\section{Proof of Theorem~\ref{th-interchange}}
\label{sec-proof-interchange}

\subsection{Basic properties of a limit of stationary distributions}

\begin{lem}
\label{lem-tightness}
The sequence of stationary distributions, i.e. the distributions of $\mathring f^n(\infty)$, is tight. Any subsequential distributional limit $\mathring f(\infty)$ is such that: \\
$\mathring f(\infty)$ is concentrated on $\mathring \cm$; 
\beql{eq-moment-bound}
\E \Phi_{1+\chi}(\mathring f(\infty)) \le \bar C,
\eeql
where $\bar C$ is the constant in Theorem~\ref{th-closeness};
\beql{eq-first-moment-bound}
\E \Phi_1(\mathring f(\infty)) \le 
%C_1 \doteq 
1+\bar C.
\eeql
\end{lem} 

{\em Proof.}  By Theorem~\ref{th-closeness} and Markov inequality,
for any $\delta>0$, there exists a sufficiently large number $C(\delta)>0$, such that for all large $n$,
with probability at least $1-\delta$, $\Phi_{1+\chi}(\mathring f^n(\infty)) \le C(\delta)$. The set
$$
S(\delta)= \{f \in \mathring \cm ~|~ \Phi_{1+\chi}(f) \le C(\delta)\} 
$$
is compact in $\cm$ (equipped with the weak convergence topology). And we know that, for all large $n$, $\pr\{\mathring f^n(\infty) \in S(\delta)\} \ge 1-\delta$. Therefore, the sequence of distributions of $\mathring f^n(\infty)$ is tight.
Consider any subsequential distributional limit $\mathring f(\infty)$, i.e. $\mathring f^n(\infty) \Rightarrow \mathring f(\infty)$
for a subsequence of $n$. 
Then
$\pr\{\mathring f(\infty) \in S(\delta)\} \ge 1-\delta$, and then the mean $\overline {\mathring f}(\infty)=0$ w.p. 1.
Moreover, the limit $\mathring f(\infty)$ must be such that \eqn{eq-moment-bound} holds.
(This is by Fatou lemma, because, using Skorohod representation,
we can construct the sequence of $\mathring f^n(\infty)$ and $\mathring f(\infty)$ on a common 
probability space such that the convergence
$\mathring f^n(\infty) \stackrel{w}{\rightarrow} \mathring f(\infty)$ is w.p.1.) Finally, \eqn{eq-first-moment-bound} is from
$$
\E \Phi_1(\mathring f(\infty)) \le \E \left[[\Phi_{1+\chi}(\mathring f(\infty))]^{1/(1+\chi)}\right] \le 1+ \E \Phi_{1+\chi}(\mathring f(\infty))
\le 1+\bar C.
$$
$\Box$

\subsection{Characterization of a limit of stationary distributions.}

Suppose $f(0) \in \mathring \cm$, that is $f(0) \in \cm$ and $\bar f(0) =0$.
If $f(\cdot)$ is a mean-field model with initial state $f(0)$ (i.e., the solution to \eqn{eq-pde-oper}),
then the corresponding centered trajectory $\mathring f(\cdot)$ will be called the {\em centered mean-field model}. 

\begin{lem}
\label{lem-stat-distr-of-mfl}
The distribution of any subsequential limit $\mathring f(\infty)$ is a stationary distribution of the 
deterministic process evolving along centered mean-field limits.
\end{lem}

{\em Proof.}  
By Theorem~\ref{th-finite-interval}(ii),
the dependence of the deterministic trajectory $\mathring f(\cdot)$ on the initial state $\mathring f(0)$ is continuous. 
Then we can apply Theorem 8.5.1 in \cite{Liptser_Shiryaev}, adapted to our setting.
Or, the proof is also easy to obtain directly as follows.
We need to show that for any test function $h \in \cc_b$ and any $t\ge 0$, we have
\beql{eq-inv-lim}
\E \mathring f(0) h = \E \mathring f(t) h,
\eeql
where $f(0)$ is equal in distribution to $f(\infty)$. We obtain \eqn{eq-inv-lim} by taking the limit of the equality
\beql{eq-inv-prelim}
\E \mathring f^n(0) h = \E \mathring f^n(t) h,
\eeql
where $\mathring f^n(0)$ is equal in distribution to $\mathring f^n(\infty)$; \eqn{eq-inv-prelim} clearly holds for all $n$.
Clearly, $\E \mathring f^n(0) h \to \E \mathring f(0) h$. To demonstrate
\beql{eq-rhs-conv}
\E \mathring f^n(t) h \to \E \mathring f(t) h,
\eeql
we can use Skorohod representation, so that the convergence
$\mathring f^n(0) \stackrel{w}{\rightarrow} \mathring f(0)$ is a.s. 
For any deterministic converging sequence $\mathring f^n(0) \stackrel{w}{\rightarrow} \mathring f(0)$ we have,
by Theorem~\ref{th-finite-interval}(ii), $\mathring f^n(t) \Rightarrow \mathring f(t)$ (which is equivalent to 
$\mathring f^n(t) \stackrel{\pr}{\rightarrow} \mathring f(t)$), and then
$\E \mathring f^n(t) h \to \E \mathring f(t) h$. Thus, we obtain 
\eqn{eq-rhs-conv}, and then \eqn{eq-inv-lim}.
$\Box$

\subsection{Completion of the proof of Theorem~\ref{th-interchange}.}

Consider any subsequential distributional limit $\mathring f(\infty)$. Its distribution is a stationary distribution of the 
deterministic process evolving along centered mean-field limits $\mathring f(\cdot)$.

{\em First proof.} Consider any two different initial conditions $\mathring f^{(1)}(0)\in \mathring \cm$ and 
$\mathring f^{(2)}(0) \in \mathring \cm$ of the deterministic process $\mathring f(\cdot)$. The Wasserstein $W_1$-distance $\|\mathring f^{(1)}(t) - \mathring f^{(2)}(t)\|_1$ must strictly decrease with $t$, by Theorem 2 in \cite{GSS96}. Then, the only option is that $\mathring f(\infty)$ is concentrated on a single element $\phi$, which then must be a traveling wave shape, and then unique traveling wave shape.
Otherwise, if we consider  two independent stationary versions of the process 
$\mathring f(\cdot)$, say $\mathring f^{(1)}(0)$ and $\mathring f^{(2)}(0)$, then $\E \|\mathring f^{(1)}(t) - \mathring f^{(2)}(t)\|_1$
is finite for $t=0$ (which follows from $\E \Phi_{1+\chi}(\mathring f(\infty)) \le \bar C$) and it is strictly decreasing in $t$;
this contradicts the stationarity.
Finally, 
we obtain \eqn{eq-tws-moment}, i.e. 
$\Phi_{1+\chi}(\phi) <\infty$, because $\E \Phi_{1+\chi}(\mathring f(\infty)) \le \bar C$. $\Box$

{\em Second proof.}  By Theorem 3.2 in \cite{St2020-wave}, we have the existence of the unique traveling wave shape 
$\phi \in \mathring \cm$. Then, by Theorem 2 in \cite{GSS96}, for any $\mathring f(0) \in \mathring \cm$, the 
Wasserstein $W_1$-distance 
$\|\mathring f(t) - \phi\|_1$ is strictly decreasing when it is non-zero. Then $\mathring f(\infty)$ must be concentrated at $\phi$.
Property \eqn{eq-tws-moment} follows from $\E \Phi_{1+\chi}(\mathring f(\infty)) \le \bar C$. $\Box$

\begin{remark}
We have given two proofs which complete the proof
of Theorem~\ref{th-interchange}. They both require additional assumptions \eqn{eq-moment-finite222}
and \eqn{eq-positive-density}. 
Note that, as a byproduct of the first proof, we also obtain the existence of the unique traveling wave shape, but only under these additional conditions. As far as the existence of the traveling wave shape is concerned, it is established in Theorem 3.1 in \cite{St2020-wave} under weaker conditions, only requiring the finite second moment of jump size, $m^{(2)} = \E Z^2 < \infty$.
\end{remark}

\section{Discussion}
\label{sec-discuss}

Main results of this paper, in a sense, complete the ``program''  represented by previous work 
\cite{GSS96,GMP97,St2020-wave} on the specific model in this paper. Paper \cite{GMP97}, informally speaking,
proves the convergence to a deterministic mean-field model as $n\to\infty.$ (In our paper we generalize that result to 
a more general model, without the finite-dependence assumption.)  Paper \cite{GSS96}, informally speaking,
proves that if a traveling wave exists, then each mean-field model trajectory is attracted to that traveling wave, as $t\to\infty$;
paper \cite{St2020-wave} shows that a traveling wave does exist under very general assumptions.
This paper proves that the convergence to the traveling wave also holds if we ``interchange the limits:'' 
first consider the stationary distribution (take limit in $t\to\infty$) and then consider the $n\to\infty$ limit 
of stationary distribution; if we take limits in this order, the limit is same -- a traveling wave. Thus, the results of this ``program,'' answer essentially ``all'' questions about the behavior of the system when $n$ is large -- both about its transient behavior 
and about its stationary distribution.

There many other well-motivated large-scale particle systems (cf. 
\cite{Man-Sch-2005, Mal-Man-2006, Manita-2009, Manita-2006, Malyshkin-2006,Manita-2014,balazs-racz-toth-2014}),
for which implementing a similar program would be of interest.

%\iffalse
%%%%%%%%%%%%\bibliographystyle{acmtrans-ims}
%%%%%%%%%%%%\bibliographystyle{apt}
\bibliographystyle{abbrv}

%\bibliography{biblio-stolyar}
%\fi

\end{document}